\documentclass[11pt,fullpage]{article}
\input epsf
\usepackage{amsfonts}
\usepackage{amscd}
\usepackage{amsmath,amssymb}
\usepackage{amstext}
\usepackage{amscd}
\usepackage[all]{xy}
\usepackage{pstricks,pst-node,pst-tree}

\newtheorem{lem}{LEMMA}[section]
\newtheorem{theo}[lem]{THEOREM}

\newtheorem{coro}[lem]{COROLLARY}
\newtheorem{prop}[lem]{PROPOSITION}

\newtheorem{definition}[lem]{Definition}
\newtheorem{rems}[lem]{Remarks}
\newtheorem{rem}[lem]{Remark}
\newtheorem{ex}[lem]{Example}
\newtheorem{exs}[lem]{Examples}
\renewcommand{\descriptionlabel}[1]%
       {\hspace{\labelsep}\textsf{#1}}

\begin{document}

\title{ Wild Knots as limit sets of Kleinian Groups}
\author{Gabriela Hinojosa\thanks{This work was partially supported by
CONACyT (Mexico), grant G36357-E, and DGEP-UNAM (M\'exico).}\\
Facultad de Ciencias, UAEM\\
62210, Cuernavaca, M\'exico.\\
gabriela@servm.fc.uaem.mx}
\date{May 17$^{th}$, 2004}

\maketitle{}
\centerline{Dedicated to Alberto Verjovsky on the occasion of his $60^{th}$ anniversary.}

\begin{abstract}
In this paper we study  kleinian
groups of Schottky type whose limit set is a wild knot in the sense of Artin
and Fox. We show that, if 
the ``original knot'' fibers over the circle then the wild knot
$\Lambda$ also
fibers over the circle. As a consequence, the universal covering of
$\mathbb{S}^{3}-\Lambda$ is $\mathbb{R}^{3}$. We prove that the
complement of a dynamically-defined fibered wild knot can not be a
complete hyperbolic 3-manifold.    
\end{abstract}

\section{Introduction}
Kleinian and (quasi)-fucshian groups were originally defined by H. Poincar\'e
in the 1880's.
In particular, quasi-fucshian groups acting on the Riemann sphere
provide us with a most natural way to construct fractal circles with the
property of being self-similar. Kleinian groups in higher dimensions
may be defined as subgroups of M\"oebius groups acting conformally on the sphere
$\mathbb{S}^{n}$. B. Maskit (see \cite{maskit}), 
M. Kapovich (see \cite{kap1}), have described kleinian groups of
schottky type whose limit set is a wild knot in the sense of Artin and
Fox.\\

The world of wild topology born with the work of
Antoine and Alexander in the 1920's, who discovered that a Cantor
set and 2-sphere (respectively) can be knotted in $\mathbb{S}^{3}$. Many works in this
direction have be done since then. However, the classical knot theory
has excluded the wild case.
Recently, Montesinos (\cite{montesinos1}, \cite{montesinos2}) has
exhibited a relationship between open 3-manifolds and wild knots and
arcs. This generalizes the relationship between closed manifolds and
tame links. Montesinos also proved that there exists a universal wild
knot, i.e. every closed orientable 3-manifold is a 3-fold
branched covering of $\mathbb{S}^{3}$ with branched set a wild
knot. This shows how rich the wild knot theory can be.\\

The purpose of this paper is to study the topological properties of
dynamically-defined wild knots. 
In the first sections we describe geometrically the 
action of these kleinian groups, their limit sets and fundamental domains.
In section 4, we study the complement $\Omega$, of dynamically-defined wild knots.
We prove that if the {\it template} (see definition 2.2) of the
corresponding wild knot is a
non-trivial tame fibered knot, then the complement of it also
fibers over the circle with fiber an infinite-genus surface with 
one end. As a consequence, we obtain that the universal covering of $\Omega$
is $\mathbb{R}^{3}$. This is not a trivial fact, e.g. there are Whitehead manifolds which are infinite
cyclic coverings spaces of other non-compact 3-manifolds (see, for
instance \cite{myers}). Although 
a generalization of a classical theorem of R. H. Bing characterizes $\mathbb{R}^{3}$ among all 
contractible open 3-manifolds (\cite{doyle}), its application to prove that $\mathbb{R}^{3}$ is the 
universal covering of the complement of any knot is not entirely
obvious, at least for us. 
As far as we know, this question remains open. 
P. H. Doyle
(\cite{doyle1}) gave an example of a wild knot whose fundamental group is $\mathbb{Z}*G$, where 
$G\neq\mathbb{Z}$ is a group. We suspect that the universal covering of
its complement is a contractible space not
homeomorphic to $\mathbb{R}^{3}$.
In section 5, we describe completely the monodromy of wild knots 
whose complements fiber over the circle. The monodromy allows us to
recognize if two dynamically-defined wild knots are isotopic.\\

In section 6, we consider the question of whether the complement of a
dynamically-defined, fibered wild knot 
can be given the structure of a complete hyperbolic 3-manifold. We recall that W. Thurston
(\cite{Thurston2}) proved that if $K$ is a tame knot, then
$\mathbb{S}^{3}-K$ admits a complete hyperbolic structure if and only if $K$ is
neither a torus knot nor a satellite knot. This motivates us to discover if
there exists an equivalent result for dynamically-defined wild knots. In this direction, we
answer the following question formulated by E. Ghys: Let $K$ be a tame knot such that its
complement is a complete hyperbolic 3-manifold. Is the complement of the wild knot obtained from $K$ still 
a complete hyperbolic 3-manifold?\\
 
In the last section, we give the definition of $q$-fold cyclic covering space
over $\mathbb{S}^{3}$ with branched set a dynamically-defined wild knot $\Lambda$, when
the ``original knot'' fibers over the circle. Recently, Montesinos 
(\cite{montesinos1}, \cite{montesinos2}) has proved amazing results 
about branched coverings of $\mathbb{S}^{3}$ with branched sets that
are wild knots. However, in our case the 
covering spaces are not manifolds.\\

This paper contains some results from my Ph.D. thesis which was
directed by A. Verjovsky. I would like to thank him for many councils
and pivotal suggestions. I would also like to thank M.
Kapovich for his very helpful corrections and comments, F. Gonz\'alez-Acu\~na and
V. N\'u\~nez for the valuable discussions.\\ 
  
\section{Preliminaries}

Let $M\ddot{o}b(\mathbb{S}^{n})$ denote the group of M\"obius
transformations of the n-sphere 
$\mathbb{S}^{n}=\mathbb{R}^{n}\cup\{\infty\}$. For a discrete group 
$G\subset M\ddot{o}b(\mathbb{S}^{n})$ {\it the discontinuity set} is
defined 
$$
\Omega (G)=\{x\in\mathbb{S}^{n}: 
x\hspace{.2cm} \mbox{has a neighbourhood} 
\hspace {.2cm}U(x)\hspace{.2cm}\mbox{such that}
$$
$$
\hspace{1cm}
U(x)\cap g(U(x))\hspace{.2cm} 
\mbox{is empty for all but finite elements}\hspace{.2cm} g\in G\}
$$
The complement 
$\mathbb{S}^{n}-\Omega(G)=\Lambda(G)$ is the {\it  limit set} (see 
\cite{kap1}). A subgroup $G\subset M\ddot{o}b(\mathbb{S}^{n})$ is
called a {\it  kleinian group} if 
$\Omega (G)$ is not empty. We will be concerned with very specific
kleinian groups of schottky type (see \cite{maskit}, page 82).\\

We recall that a conformal map $\psi$ on $\mathbb{S}^{n}$ can be extended in a
natural way to the hyperbolic space  $\mathbb{H}^{n+1}$ 
as an orientation-preserving isometry with
respect to the Poincar\'e metric. Hence we can identify 
the group $M\ddot{o}b(\mathbb{S}^{n})$ with the group of orientation-preserving isometries of the
hyperbolic $(n+1)$-space $\mathbb{H}^{n+1}$. Thus
we can also define the limit set of a kleinian group through sequences (see \cite{maskit} section II.D).\\ 
 
\begin{definition}
A point $x$ is a {\it  limit point} for the kleinian group
$G$ if there exist a point $z\in\mathbb{S}^{n}$ and a sequence $\{g_{m}\}$ 
of {\it distinct elements} of $G$, with $g_{m}(z)\rightarrow x$. 
The set of limit points is $\Lambda (G)$.
\end{definition}

One way to illustrate the action of a kleinian group $G$ is to draw a picture 
of $\Omega/G$. For this purpose a fundamental domain is very helpful. Roughly 
speaking, it contains one point from each equivalence class in $\Omega$ (see 
\cite{kap2} pages  78-79, \cite{maskit} pages  29-30). Let $D$ be the fundamental
  domain of $G$ and consider the orbit space  
$D^{*}=\overline{D}\cap\Omega/\sim_{G}$ with
the quotient topology. Then $D^{*}$ is homeomorphic to $\Omega/G$.

\begin{definition}
A  necklace $T$ of $n$-pearls ($n\geq 3$), is a collection of $n$ consecutive 
2-spheres $\Sigma_{1}$, $\Sigma_{2},\ldots ,\Sigma_{n}$ in
$\mathbb{S}^{3}$, such that $\Sigma_{i}\cap\Sigma_{j}=\emptyset$
($i\neq i+1, i-1$ mod $n$), except that
$\Sigma_{i}$ and $\Sigma_{i+1}$ are tangent $i=1,2,\ldots,n-1$ and
$\Sigma_{1}$ and $\Sigma_{n}$ are 
tangent. Each
2-sphere is called a {\it  pearl}. 
\end{definition}
\begin{rem}
If the points of tangency are
joined by spherical geodesic segments in $\mathbb{S}^{3}$, we obtain a 
tame knot $K$ which is called the {\it template} of $T$. Observe that 
the ordering of pearls gives us an orientation for $K$.\\

Conversely, a pearl-necklace $T$ subordinate to a polygonal knot $K$, is a collection of round 
2-spheres $\Sigma_{1}$, $\Sigma_{2},\ldots ,\Sigma_{n}$ in
$\mathbb{S}^{3}$, such that $\Sigma_{i}$ is tangent to 
$\Sigma_{i+1}$ (index $i$ mod $n$) and  $K$ is totally covered by them. For simplicity, we also require that
each segment of $K$ contained in the interior of each pearl is an unknotted tangle.
\end{rem}

We define the {\it filling of $T$} as $|T|=\cup^{n}_{i=1} B_{i}$, 
where $B_{i}$ is the round closed 3-ball whose boundary $\partial
B_{i}$ is $\Sigma_{i}$.

\begin{ex}
 $K=$ Trefoil knot.
\begin{figure}[tbh]
\centerline{\epsfxsize=1.3in \epsfbox{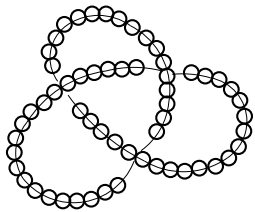}}
\caption{\sl A pearl-necklace whose template is the trefoil knot.}
\label{F1}
\end{figure}
\end{ex}

Let $\Gamma$ be the group generated by reflections $I_{j}$, through 
$\Sigma_{j}$ ($j=1,\ldots,n$). Then $\Gamma$ is a kleinian group
of Schottky type.
It is easy to verify that the fundamental domain for the kleinian group $\Gamma$ is 
$D=\mathbb{S}^{3}-|T|$.

\section{Geometric Description of the Limit Set}

To describe the limit set of $\Gamma$, we need to find all the
 accumulation points of its  orbits. We
shall thus consider all the possible sequences of elements of   
$\Gamma$ . We will do this in stages:

\begin{enumerate}

\item First stage: Observe that each reflection map $I_{j}$ 
($j=1,2,\ldots,n$), sends a copy of the exterior of $K$ into the ball $B_{j}$. In particular,
a mirror image of $K$ is sent into $B_{j}$ (see Figure 2). 

\begin{figure}[tbh]
\centerline{\epsfxsize=1.3in \epsfbox{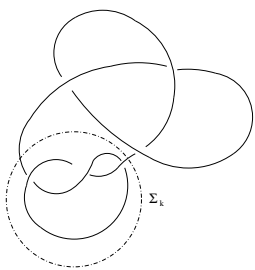}}
\caption{\sl Image of the reflection map $I_{k}$.}
\label{F2}
\end{figure}

After doing this for each $j$, we obtain a new necklace
$T_{1}$ of $n(n-1)$ pearls, subordinate to a new knot $K_{1}$; which
is in turn isotopic to the connected sum of $n$ copies of the mirror image of $K$ and
one copy of $K$.

\item Second stage: Now, we reflect with respect to each pearl of
$T_{1}$. When we are finished, we obtain a new necklace $T_{2}$ of $n(n-1)^{2}$ pearls, whose 
template is a tame knot $K_{2}$; which is in turn isotopic to the connected
sum of $n^{2}-n+1$ copies of $K$ and $n$ copies of its mirror image (recall that composition
of an even number of reflections is orientation-preserving). Observe that  
$|T_{2}|\subset |T_{1}|$ (see Figure 3).

\begin{figure}[tbh]
\centerline{\epsfxsize=1.3in \epsfbox{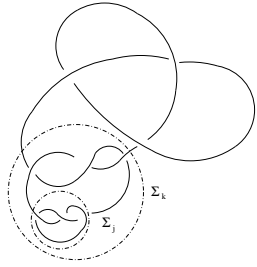}}
\caption{\sl Image of the reflection map $I_{k}$
         after reflecting with respect to $\Sigma_{j}$.}
\label{F4}
\end{figure}

\item $k^{th}$-stage: In this case, we reflect with respect to each pearl of $T_{k-1}$, and we obtain 
a necklace $T_{k}$ of $n(n-1)^{k}$ pearls, subordinate to a tame 
knot $K_{k}$. By construction, $|T_{k}|\subset |T_{k-1}|$.
\end{enumerate}

Let $x\in\cap_{k=1}^{\infty} |T_{k}|$. We shall prove that $x$ is a
 limit point. Indeed, there exists a sequence of closed balls
 $\{B_{m}\}$ with $B_{m}\subset |T_{m}|$ such that $x\in B_{m}$ for
 each $m$. For any $z\in\mathbb{S}^{3}-|T|$ we can find
 a sequence $\{w_{m}\}$ of distinct elements of $\Gamma$, such
 that $w_{m}(z)\in B_{m}$. Since $diam(B_{m})\rightarrow 0$, it follows
 that $w_{m}(z)$ converges to $x$ independently of $z$. Then, $x$ is a limit point of $\Gamma$ (see
definition 2.1). 
The other inclusion clearly holds, since the fundamental domain of $\Gamma$ is
$\mathbb{S}^{3}-|T|$. Therefore, the limit set is given by
$$
\Lambda(\Gamma)=\varprojlim_{k} |T_{k}|=\bigcap_{k=1}^{\infty} |T_{k}|.
$$

\begin{lem} (\cite{maskit}). Let $T$ be a pearl-necklace subordinate
to the non-trivial, tame knot $K$. Let $\Gamma$ be the group generated by reflections through each pearl.
Then the limit set $\Lambda(\Gamma)$ is homeomorphic to  $\mathbb{S}^{1}$. 
\end{lem} 

\begin{lem} Let $T$ be a pearl-necklace
 subordinate to the non-trivial tame knot $K$. Then
 $\Lambda(\Gamma)$ is wildly embedded in
 $\mathbb{S}^{3}$. Moreover, $\Lambda(\Gamma)$ is wild at every point.
\end{lem}

\hspace{-.67cm}{\bf {\it Proof.}} (Compare \cite{maskit}). Let $x\in
\Lambda(\Gamma)$ be a limit point. Since the
knot $K$ keeps reproducing either itself or its mirror image, we have that given an open,
connected neighborhood $U$ of $X$, there are infinitely many copies $K$ and its mirror image 
in $U$.\\

We will denote the complement of the wild knot $\Lambda$ contained in $U$ by
$U-\Lambda(\Gamma)$, and the complement of the pearl-necklace $T_{k}$ in $U$ by $U-|T_{k}|$.\\
 
Hence, using Van-Kampen's Theorem and the fact that
$\Pi_{1}(U-\Lambda(\Gamma))$ 
is the direct limit of 
$\{\Pi_{1}(U-|T_{k}|),\hspace{.2cm}
k=l,l+1,\ldots;\hspace{.3cm}f_{k},\hspace{.2cm}k=l,l+1,\ldots)\}$.
Where
$f_{k}:\Pi_{1}(\mathbb{S}^{3}-|T_{k}|)\rightarrow\Pi_{1}(\mathbb{S}^{3}-|T_{k+1}|)$
is the inclusion map. The index $l$ is the smallest integer such that
some balls of $T_{l}$ are totally contained in $U$ (see Lemma 2.4.1 in \cite{rushing}). We have that,
$$
\Pi_{1}(U-\Lambda(\Gamma))
\cong(\cdots((\Pi_{1}(K)*_{\{z\}}\Pi_{1}(K))*_{\{z\}}\cdots *_{\{z\}}\Pi_{1}(K))
*_{\{z\}} \cdots,
$$
is an infinite free product of the fundamental group of the knot with itself. This implies that
$\Pi_{1}(U-\Lambda(\Gamma))$  is not isomorphic 
to a finitely generated group, i.e. it
is infinitely generated. Therefore, the result follows. $\blacksquare$

\section{Fibration of $\mathbb{S}^{3}-\Lambda(\Gamma)$ over
 $\mathbb{S}^{1}$}

We recall that a knot or link $L$ in $\mathbb{S}^{3}$ is {\it fibered} if there exists a 
locally trivial fibration $f:(\mathbb{S}^{3}-L)\rightarrow \mathbb{S}^{1}$. We 
require that $f$ be well-behaved near $L$, that is, each component $L_{i}$ is 
to have a neighbourhood framed as $\mathbb{D}^{2}\times\mathbb{S}^{1}$, with 
$L_{i}\cong \{0\}\times\mathbb{S}^{1}$, in such a way that the
restriction of $f$ to $(\mathbb{D}^{2}-\{0\})\times\mathbb{S}^{1}$ is the map 
into $\mathbb{S}^{1}$ given by $(x,y)\rightarrow \frac{y}{|y|}$. It follows that each
$f^{-1}(x)\cup L$, $x\in\mathbb{S}^{1}$, is a 2-manifold
 with boundary $L$: in fact a Seifert surface for $L$ 
(see \cite{rolfsen}, page 323).

\begin{exs}
The right-handed trefoil knot and the figure-eight knot are fibered knots with fiber
the punctured torus.
\end{exs}

\begin{lem}
Let $T$ be an $n$-pearl necklace subordinate to the fibered tame knot $K$,
with fiber $S$. Then
$\Omega (\Gamma)/\Gamma$ fibers over the circle with fiber $S^{*}$,
the closure of $S$ in $\mathbb{S}^{3}$ with $n$ boundary points removed.
\end{lem}

\hspace{-.67cm}{\bf {\it Proof.}}
Let $\widetilde{P}:(\mathbb{S}^{3}-K)\rightarrow \mathbb{S}^{1}$ be the
given fibration with fiber the surface $S$. Observe that 
$\widetilde{P}\mid_{\mathbb{S}^{3}-|T|}\equiv P$ is a fibration and, after
modifying $\widetilde{P}$ by isotopy if necessary, we can consider  that
the fiber $S$  cuts each pearl $\Sigma_{i}\in T$ in arcs
$a_{i}$, whose end-points are $\Sigma_{i-1}\cap\Sigma_{i}$ and $\Sigma_{i}\cap\Sigma_{i+1}$. These two
points belong to the limit set (see Figure 4).

\begin{figure}[tbh]
\centerline{\epsfxsize=1.2in \epsfbox{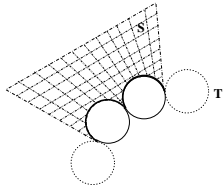}}
\caption{\sl The fiber intersects each pearl in arcs.}
\label{F6}
\end{figure}

Hence the space $\overline{D}\cap\Omega(\Gamma)$ fibers over the
circle with fiber the 2-manifold $S^{*}$, which is the closure of the
surface $S$ in $\mathbb{S}^{3}$ with $n$ boundary points removed
(points of tangency of the pearls). Since  
$\Omega(\Gamma)/\Gamma\cong D^{*}$ and in our case 
$D^{*}=\overline{D}\cap\Omega(\Gamma)$, the 
result follows. $\blacksquare$\\

By the above Lemma, in order to describe completely the orbit space 
$(\mathbb{S}^{3}-\Lambda(\Gamma))/\Gamma$ 
in the case that the original knot 
is fibered, we only need to determine its monodromy: which is
precisely the monodromy of the knot.\\

Consider the orientation-preserving index two subgroup $\widetilde{\Gamma}\subset\Gamma$. Since 
$\widetilde{\Gamma}$ is a normal subgroup of $\Gamma$, it follows by Lemma 8.1.3 in 
\cite{thurston1} that $\widetilde{\Gamma}$ has the same limit set that $\Gamma$. Therefore
$\mathbb{S}^{3}-\Lambda(\Gamma)=\mathbb{S}^{3}-\Lambda(\widetilde{\Gamma})$. 

\begin{lem}
Let $T$ be an $n$-pearl necklace subordinate to the fibered tame knot $K$
with fiber $S$. Then
$\Omega (\widetilde{\Gamma})/\widetilde{\Gamma}$ fibers over the circle with fiber a
surface $S^{**}$, which is homeomorphic to the surface
$S^{*}$ joined along an arc to a copy of itself in $\mathbb{S}^{3}$,
with $2(n-1)$ boundary points removed
(the points of tangency).
\end{lem}

\hspace{-.67cm}{\bf {\it Proof.}}
It is easy to see that the fundamental domain for $\widetilde{\Gamma}$ is 
$$
\widetilde{D}=(\mathbb{S}^{3}-|T|)\cup (B_{j}-I_{j}(|T-\Sigma_{j}|)),
$$
for some $j$. The rest of the proof is very similar to the proof of the above
Lemma. $\blacksquare$\\

\begin{theo}
Let $\Sigma_{1}$, $\Sigma_{2},\ldots ,\Sigma_{n}$ be round 2-spheres
in $\mathbb{S}^{3}$ which form a necklace whose template is a non-trivial
tame fibered knot. Let $\Gamma$ be the group generated by reflections
$I_{j}$ on $\Sigma_{j}$ ($j=1,2,\ldots,n$) and let $\widetilde{\Gamma}$
be the orientation-preserving index two subgroup of $\Gamma$. 
Let $\Lambda(\Gamma)=\Lambda(\widetilde{\Gamma})$ be the 
corresponding limit set. Then:
\begin{enumerate}
\item There exists a locally trivial fibration $\psi
:(\mathbb{S}^{3}-\Lambda(\Gamma))\rightarrow\mathbb{S}^{1}$, where the
 fiber $\Sigma_{\theta}=\psi^{-1}(\theta)$ is an orientable
 infinite-genus surface with one end.
\item  $\overline{\Sigma_{\theta}}-\Sigma_{\theta}=\Lambda(\Gamma)$, where
$\overline{\Sigma_{\theta}}$ is the closure of $\Sigma_{\theta}$ in $\mathbb{S}^{3}$.
\end{enumerate}
\end{theo}

\hspace{-.67cm}{\bf {\it Proof.}}
We will first prove that $\mathbb{S}^{3}-\Lambda(\Gamma)$ fibers over the circle.
We know that $\zeta:\Omega(\widetilde{\Gamma)}\rightarrow\Omega(\widetilde{\Gamma)}/\widetilde{\Gamma}$ 
is an infinite-fold covering since $\widetilde{\Gamma}$ acts freely on $\Omega(\widetilde{\Gamma})$. By
the previous lemma, there exists a locally trivial 
fibration $\phi:\Omega(\widetilde{\Gamma})/\widetilde{\Gamma}\rightarrow\mathbb{S}^{1}$ with fiber 
$S^{**}$.
Then $\psi=\phi\circ\zeta:\Omega(\widetilde{\Gamma})\rightarrow\mathbb{S}^{1}$ is a 
locally trivial fibration. The fiber is $\Gamma(S^{**})$, i.e. the orbit of 
the fiber.\\

 We now describe $\Sigma_{\theta}=\widetilde{\Gamma}(S^{**})$ in detail.
Let $\widetilde{P}:(\mathbb{S}^{3}-K)\rightarrow \mathbb{S}^{1}$ be the
given fibration. The fibration $\widetilde{P}\mid_{\mathbb{S}^{3}-|T|}\equiv
P$ has been chosen as in the lemmas above. The 
fiber $\widetilde{P}^{-1}(\theta)=P^{-1}(\theta)$ is a Seifert surface $S^{*}$ of $K$, for each 
$\theta\in\mathbb{S}^{1}$. We suppose $S^{*}$ is oriented. Recall that the boundary of $S^{*}$ cuts 
each pearl $\Sigma_{j}$ in an arc $a_{j}$ going from one point of
tangency to another.\\

The reflection  $I_{j}$ maps both a copy of $T-\Sigma_{j}$
(called $T^{j}$) and a copy of $S^{*}$ (called $S_{1}^{*j}$) into
the ball $B_{j}$, for $j=1,2,\ldots,n$. Observe that both $T^{j}$ and $S^{*j}$ have opposite orientation 
and that $S^{*}$ and $S_{1}^{*j}$ are joined by
the arc $a_{j}$ (see Figure 5) which, in both surfaces, has the same orientation. 

\begin{figure}[tbh]
\centerline{\epsfxsize=.8in \epsfbox{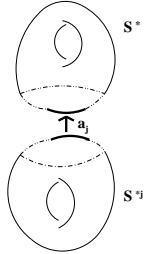}}
\caption{\sl Sum of two surfaces $S^{*}$ and $S^{*j}$ along arc $a_{j}$.}
\label{F7}
\end{figure}

The necklaces $T^{j}$ and $T$ are joined by the points of tangency of the pearl $\Sigma_{j}$ and the 
orientation of these two points is preserved by the reflection $I_{j}$. Thus, we have obtained a new 
pearl-necklace isotopic to the 
connected sum $T\#T^{j}$, whose complement also fibers over the circle with 
fiber the sum of $S^{*}$ with $S_{1}^{*j}$ along arc $a_{j}$, namely
the fiber is $S^{*}\#_{a_{j}}S_{1}^{*j}$.\\

Now do this for each $j=1,\cdots,n$. At the end of the first stage, we
have a new pearl-necklace $T_{1}$ whose template is the knot $K_{1}$
(see section 3). Its complement fibers 
over the circle with fiber the Seifert surface $S^{*}_{1}$, which is in turn 
homeomorphic to the sum of $n+1$ copies of $S^{*}$ along the respective arcs.\\

Continuing this process, we have from the second step onwards, that $n-1$ copies of $S^{*}$ are added 
along arcs to 
each surface $S^{*i}_{k}$, (the $i^{th}$ copy of $S^{*}$ corresponding to the $k^{th}$ stage). Notice that 
in each step,
the points of tangency are removed since they belong to the limit set,
and the length of the arcs $a_{j}$ tends to zero.\\

From the remarks above, we have that $\Sigma_{\theta}$ is homeomorphic to an orientable
infinite-genus surface. In fact, it is the sum along arcs of an infinite
number of copies of $S^{*}$. To determine what kind of surface it is,
according to the classification theorem of non-compact surfaces (see
\cite{richards}), we need only describe its set of ends.
\vskip .3cm
Consider the fuchsian model (see \cite{maskit}). In this case, the
necklace is formed by pearls of the same size, and each 
pearl is orthogonal to the unit circle (its template). Then its limit set
is the circle and its complement fibers over $\mathbb{S}^{1}$
with fiber the disk.
\vskip .3cm
In each step we are adding handles to this disk in such a way that
they accumulate on the boundary. If we intersect this disk with any
compact set, we have just one connected component. Hence it has only
one end. Therefore, our surface has one end.

\begin{figure}[tbh]
\centerline{\epsfxsize=1in \epsfbox{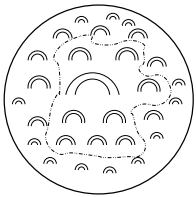}}
\caption{\sl Disk with handles intersected with a compact set.}
\label{F8}
\end{figure}

The first part of the theorem has been proved. For the second part,
observe that the closure of the fiber in $\mathbb{S}^{3}$ is the fiber
union its boundary. Therefore 
$\overline{\Sigma_{\theta}}-\Sigma_{\theta}=\Lambda(\Gamma)$. $\blacksquare$

\begin{rems}
\begin{enumerate}
\item Theorem 4.4 allows us to compute the monodromy of $\Lambda$ (see next section), and describe
completely its complement. We would like to point out that this fact provides a way
to determine if two dynamically-defined, fibered wild knots are isotopic.
\item This theorem can be generalized to  fibered links.
\item This theorem gives an open book decomposition (\cite{rolfsen}
pages  340-341) of
$\mathbb{S}^{3}-\Lambda(\Gamma)$, where the ``binding'' is the wild
knot $\Lambda(\Gamma)$, and each ``page'' is an
orientable. infinite-genus surface with one end (the fiber).\\
Indeed, this decomposition can be thought of in the following way. For the above 
theorem, $\mathbb{S}^{3}-\Lambda(\Gamma)$ is $\Sigma_{\theta}\times [0,1]$ modulo the 
identification of the top with the bottom through a characteristic 
homeomorphism. Consider $\overline{\Sigma_{\theta}}\times [0,1]$ and identify
the top with the bottom. This is equivalent to keeping  
$\partial\overline{\Sigma_{\theta}}$ fixed 
and spinning $\Sigma_{\theta}\times \{0\}$ with respect to
$\partial\overline{\Sigma_{\theta}}$ until it is
glued  with $\Sigma_{\theta}\times \{1\}$. 
Removing $\partial\overline{\Sigma_{\theta}}$ we obtain the open book decomposition. 
\item The fibering map $\psi$ can be chosen to be smooth.\\
In fact, by a theorem of Nielsen (\cite{nielsen}) we know that every surface bundle over
$\mathbb{S}^{1}$ depends only on the monodromy isotopy class and 
by another theorem of Nielsen we have that every homeomorphism of a surface is 
isotopic to a diffeomorphism.
\end{enumerate}
\end{rems}

\begin{coro} The universal covering of
$\mathbb{S}^{3}-\Lambda(\Gamma)$ is $\mathbb{R}^{3}$.
\end{coro}

{\bf {\it Proof.}}
Let  
$P:(\mathbb{S}^{3}-\Lambda(\Gamma))\rightarrow\mathbb{S}^{1}$ be a
locally trivial fibration given by the
above theorem. If $\exp:\mathbb{R}^{1}\rightarrow\mathbb{S}^{1}$
denotes the covering map $r\mapsto e^{ir}$, then the following
diagram commutes
$$
\begin{CD}
\widetilde{X}@>\widetilde{P}>>\mathbb{R}^{1}\\
@V\pi VV@V\exp VV\\
\mathbb{S}^{3}-\Lambda(\Gamma)@>P>>\mathbb{S}^{1}
\end{CD}
$$
where $\widetilde{P}:\widetilde{X}\rightarrow\mathbb{R}^{1}$ is the pull-back of the
fibration $P:(\mathbb{S}^{3}-\Lambda(\Gamma))\rightarrow\mathbb{S}^{1}$.\\

Any fibration with contractible base space is trivial, i.e. the total 
space is homeomorphic to the product of the fiber with the base space, 
$\widetilde{X}\cong\mathbb{R}^{1}\times\Sigma_{\theta}$, where 
$\Sigma_{\theta}=\mbox{fiber}$.
The universal covering of $\mathbb{S}^{1}$ is $\mathbb{R}$ and the
universal covering of $\Sigma_{\theta}$ is
$\mathbb{R}^{2}$. Therefore, the result follows. $\blacksquare$

\begin{rem}
We already know  by the Sphere Theorem (see \cite{rolfsen} page 102) that
the universal covering of 
$\mathbb{S}^{3}-\Lambda(\Gamma)$ must be contractible. 
We point out that there are uncountable contractible open 3-manifolds which are not
homeomorphic to $\mathbb{R}^{3}$ (see \cite{mcmillan}). P. H. Doyle
(\cite{doyle1}) gave an example of a wild knot whose fundamental group is $\mathbb{Z}*G$, where 
$G\neq\mathbb{Z}$ is a group. We suspect that the universal covering
of its complement is a contractible space
not homeomorphic to $\mathbb{R}^{3}$.
\end{rem}

\begin{rems}
\begin{enumerate}
\item Through the reflecting process, we have constructed a ``Seifert
Surface'' for the wild knot $\Lambda(\Gamma)$ with the property that
its interior is a regular surface and its closure is a ``crumpled
surface'' which is not, technically, a surface. This remains true for any 
template and not just for the fibered ones. Thus the wild knots obtained
by our dynamic construction all admit Seifert surfaces in this 
{\rm generalized sense}. 
\item Let
$E=\overline{\Sigma_{\theta}}\cup\overline{\Sigma_{\theta+\pi}}$. Then
$E$ is a wild equator of $\mathbb{S}^{3}$, i.e. there exists a homeomorphism 
$h:\mathbb{S}^{3}\rightarrow\mathbb{S}^{3}$ such that $h(E)=E$ and
$h(\mathcal{H}^{+})=\mathcal{H}^{-}$, where 
$\mathcal{H}^{+}$=$\{x\in\mathbb{S}^{3}-E:x\in\Sigma_{\psi},\theta<\psi<\theta+\pi\}$
and $\mathcal{H}^{-}$= $\mathbb{S}^{3}-E-\mathcal{H}^{+}$.
\end{enumerate}
\end{rems}

\section{Monodromy}

Let $K$ be a non-trivial fibered tame knot and let $S$ be the oriented fiber. Since 
$\mathbb{S}^{3}-K$ fibers over the circle, we know that
$\mathbb{S}^{3}-K$ is a mapping torus equal to $S\times [0,1]$ modulo  
a characteristic homeomorphism $\psi:S\rightarrow S$ that glues $S\times\{0\}$ with
$S\times\{1\}$. This homeomorphism induces a homomorphism
$\psi_{\#}:\Pi_{1}(S)\rightarrow\Pi_{1}(S)$
called {\it the monodromy of the fibration}.\\

Another way to understand the monodromy is through  the 
{\it Poincar\'e first return map} of a flow defined as follows. Let $M$ be
a connected, compact  
manifold and let $f_{t}$ be a flow that possesses a transversal 
section $\eta$. It follows that if $x\in\eta$ then there exists a
continuous function $t(x)>0$ such 
that $f_{t}\in\eta$. We may define the first return Poincar\'e map  
$F:\eta\rightarrow\eta$ as $F(x)=f_{t(x)}(x)$. This map is a diffeomorphism 
and induces a homomorphism of $\Pi_{1}$, 
called {\it the monodromy} (see \cite{verjovsky}, chapter 5).\\

For the manifold $\mathbb{S}^{3}-K$, the flow that defines the first return 
Poincar\'e map $\Phi$, is the flow that cuts transversally each page of its 
open book decomposition. Notice that we can choose this flow, up to
isotopy, in such a way that the first return Poincar\'e map can be
extended to $K$ as the identity.\\

Consider a pearl-necklace $T$ subordinate to $K$. As we have observed
during the reflecting process, $K$ and $S$ are copied in each
reflection (preserving or reversing orientation). So the flow $\Phi$ is also copied, and its
direction changes according to the number of reflections. Hence the 
Poincar\'e map can be extended at each stage, giving us in the limit a
homeomorphism $\psi:\Sigma_{\theta}\rightarrow \Sigma_{\theta}$ that identifies 
$\Sigma_{\theta}\times\{0\}$ with $\Sigma_{\theta}\times\{1\}$ and which induces
the monodromy of the wild knot.\\

From the above, if we know the monodromy of the knot $K$ then we know
the monodromy of the wild knot $\Lambda(\Gamma)$.\\

By the long exact sequence associated to a fibration, we have
\begin{equation}
0\rightarrow \Pi_{1}(\Sigma_{\theta})\rightarrow
\Pi_{1}(\mathbb{S}^{3}-\Lambda(\Gamma))\stackrel{\Psi}{\overset{\longleftarrow}{\rightarrow}}
\mathbb{Z}\rightarrow 0, \tag{1}
\end{equation}
which has a homomorphism section $\Psi:\mathbb{Z}\rightarrow\mathbb{S}^{3}-\Lambda(\Gamma)$. 
Therefore (1) splits. As a consequence $\Pi_{1}(\mathbb{S}^{3}-\Lambda(\Gamma))$ is
the semi-direct product of $\mathbb{Z}$ with $\Pi_{1}(\Sigma_{\theta})$.\\

Observe that by applying Van-Kampen's Theorem in each step coupled with lemma 2.4.1 in \cite{rushing},  
$\Pi_{1}(\Sigma_{\theta})$ is isomorphic to the infinite free product of $\Pi_{1}(S)$ with itself.\\

\begin{ex}
Let $K$ be the trefoil knot. Then $S$ is the punctured torus and its
fundamental group is the free group in two generators, $a$ and
$b$. The monodromy map $\psi_{\#}$ sends $a\mapsto b^{-1}$ and
$b\mapsto ab$. Its order is six up to an outer automorphism 
(See \cite{rolfsen} pages 330-333).\\

The monodromy in the limit $\psi_{\#}:\Pi_{1}(\Sigma_{\theta})\rightarrow 
\Pi_{1}(\Sigma_{\theta})$ is given by $a_{i}\mapsto b^{-1}_{i}$ and 
$b_{i}\mapsto a_{i}b_{i}$, where
$\Pi_{1}(\Sigma_{\theta})=\{a_{i},b_{i}\}$. So
$$
\begin{aligned}
\Pi_{1}(\mathbb{S}^{3}-\Lambda (\Gamma))&\cong
\Pi_{1}(\mathbb{S}^{1})\ltimes_{\psi_{\#}}\Pi_{1}(\Sigma_{\theta})\\
&=\{a_{i},b_{i},c:a_{i}*c=b_{i}^{-1},\hspace{.2cm}b_{i}*c=a_{i}b_{i}\}\\
&=\{a_{i},c:c^{-1}a_{i}^{-1}c=a_{i}c^{-1}a_{i}^{-1}\}\\
&=\{a_{i},c:c=a_{i}ca_{i}c^{-1}a_{i}^{-1}\}\\
&=\{a_{i},c:c=ca_{i}ca_{i}c^{-1}a_{i}^{-1}c^{-1}\}\\
&=\{a_{i},c:c=ca_{i}c^{-1}c^{2}a_{i}c^{-1}c^{-1}ca_{i}^{-1}c^{-1}\}.\\
\end{aligned}
$$
Let $\alpha_{i}=ca_{i}c^{-1}$;\\

$$
\begin{aligned}
\hspace{.5cm}&=\{\alpha_{i},c:c=\alpha_{i}c\alpha_{i}c^{-1}\alpha_{i}^{-1}\}\\
&=\{\alpha_{i},c:c\alpha_{i}c=\alpha_{i}c\alpha_{i}\}
\end{aligned}
$$
This gives another method for computing the fundamental group of a wild knot
whose complement fibers over the circle.
\end{ex}

\begin{coro}Let $T$ be a pearl-necklace whose template is a
non-trivial tame fibered knot, with fiber $S$. Then
$\Pi_{1}(\Omega(\Gamma)/\Gamma)\cong\mathbb{Z}\ltimes\Pi_{1}(S^{\*})$ and 
$\Pi_{n}(\Omega(\Gamma)/\Gamma)=0$ for $n>1$.
\end{coro}

\section{Hyperbolicity}

In this section, we consider the question of whether  the complement of a dynamically-defined,
fibered wild knot can be given the structure of a complete 
hyperbolic 3-manifold.

\begin{theo}
Let $K$ be a fibered tame knot, whose complement admits a complete
hyperbolic structure. Then the complement of the 
wild knot $\Lambda(\Gamma)$, obtained from $K$ through a reflecting process is not hyperbolic.
\end{theo}

\hspace{-.67cm}{\bf {\it Proof.}}
Suppose that $\Omega=\mathbb{S}^{3}-\Lambda$ is a complete hyperbolic
3-manifold. By virtue of the way 
$\Lambda$ was constructed (see section 3), we can decompose $\Omega$ into pieces 
$W_{k}=\mathbb{S}^{3}-|T_{k}|$, such that $\Omega=\cup W_{k}$. Now,
$i:W_{2}\hookrightarrow \Omega$ 
induces an injective map  $i_{\#}:\Pi_{1}(W_{2})\hookrightarrow
\Pi_{1}(\Omega)$. If we denote the image of this group by
$G=i_{\#}(\Pi_{1}(W_{2}))$, then  $M=\mathbb{H}^{3}/G$ is a complete
hyperbolic 3-manifold which covers $W_{2}$. Let $p:M\rightarrow W_{2}$ be this covering\\

Since $G$ is finitely generated, by a theorem of Scott (see
\cite{finites}) it follows that there is a 
compact submanifold $M_{T}$ of $M$, such that the inclusion 
map $i:M_{T}\rightarrow M$ induces an isomorphism of fundamental
groups. In our case $M_{T}$ is
homeomorphic to the closure of the copy of  
$W_{2}$ in $M$. Moreover, there is a bijective
correspondence between the boundary components of $M_{T}$ and the topological ends of $M$. The 
$\mathbb{Z}\bigoplus\mathbb{Z}$-cusps correspond precisely to the toroidal components of $\partial M_{T}$.\\

By the above, $M$ has one end which is a
$\mathbb{Z}\bigoplus\mathbb{Z}$-cusp and
$M$ is diffeomorphic to $\mathbb{T}^{2}\times\mathbb{R}$. Any 2-torus  $T^{2}_{t}:=\mathbb{T}^2 \times{t}$,
$t\in\mathbb {R}$, separates $M$ into two pieces $A_{t}$ with finite volume
and  $B_{t}$. One can assume, without loss of generality,
that the hyperbolic diameter of $T^{2}_{t}$ tends to zero as $t\to\infty$. This implies that
$p(A_{t})\subset W_{2}$ has finite volume and $p(A_{t})$ determines the only end
of $W_{2}$. Therefore $W_{2}-p(A_{t})$ has finite volume and hence
$W_{2}$ itself is of finite volume (compare \cite{finites} page 253).\\

Now, let $V$ be a closed tubular
neighbourhood of $|T|\subset W_{2}$. The boundary of $V$ is an
incompressible torus in $W_{2}$ and 
it is not isotopic to a component of its 
boundary. This contradicts  Hyperbolization Theorem of Thurston. Therefore the result follows. 
$\blacksquare$\\

\begin{rem}
We recall that a complete hyperbolic manifold $N$ is said to be tame
if it is homeomorphic to the interior of a compact 3-manifold.
Marden's Tameness Conjecture says that a complete hyperbolic manifold with
finitely generated fundamental group, is tame.
Recently, Ian Agol has announced a proof of this conjecture 
(http://atlas-conferences.com/cgi-bin/abstract/camc-67).
In our case this implies immediately that $W_{2}$ has finite volume.
\end{rem}

On the other hand, we have the following positive result.
\begin{theo}
Let $K$ be a fibered tame knot whose complement is a complete hyperbolic 3-manifold. 
Then, the complement of the limit set $\Lambda(\Gamma)$ obtained from
$K$, admits a canonical decomposition into a countable union of submanifolds
$M_{j}$ with boundary. Their interiors are 
pairwise disjoint, and each interior of  $M_{j}$ is
hyperbolic. The boundary of $M_{j}$ is composed by a finite number of 
pairwise disjoint properly embedded incompressible cylinders. 
\end{theo}

\hspace{-.67cm}{\bf {\it Proof.}}
Let $T$ be a necklace of $n$ pearls subordinate to $K$. Since 
$\mathbb{S}^{3}-|T|$ is the
mapping torus of a homeomorphism $\psi: S\rightarrow S$, 
where $S$ is the fiber. Then  by the Hyperbolization Theorem  
(see \cite{mcmullen} chapter 3), it follows that $\psi$ is pseudo-Anosov.\\

As we discussed in the section 5, when we reflect with respect to 
a pearl $\Sigma_{i}$, both the monodromy $\psi$ and $T-\Sigma_{i}$ are copied into  
$Int(B_{i})$, where $B_{i}$ is the ball whose boundary is $\Sigma_{i}$. 
Since $Int(B_{i})-I(|T-\Sigma_{i}|)$ also fibers over the circle and $\psi$ is pseudo-Anosov, we have
that it is hyperbolic (see \cite{mcmullen} chapter 3). Notice that the hyperbolic pieces are
equivalent to the complement of $K$. They are separated by the corresponding pearls without
the two points of tangency, and clearly they are incompressible
cylinders. Hence the result follows. $\blacksquare$\\ 

Observe that in the first step of the above decomposition, $n$
copies of $M_{1}$ are added to it along cylinder components of its 
boundary. From the second step onwards, $n-1$ copies are added along 
cylinder components to each $M_{j}$ of the previous step.\\

In \cite{Thurston2} W. Thurston stated some open questions concerning
3-manifolds and kleinian groups. The first one is: Do
all 3-manifolds have decompositions into geometric pieces?. In this
context, as a consequence of the above Theorem, we have the following

\begin{coro}
Let $K$ be a tame, fibered knot, whose complement is a hyperbolic 3-manifold. Then
the complement of the limit set obtained from $K$, has a canonical 
decompositions into geometric pieces. In other words, it satisfies a
recasting of the  Thurston's
Geometrization Conjecture for non-compact manifolds.
\end{coro}

\section{Cyclic Coverings}

Let $K$ be a non-trivial tame fibered knot with fiber the surface
$S$. Then
$\mathbb{S}^{3}-\Lambda(\Gamma)$ is the mapping torus of a
homeomorphism $\psi:\Sigma_{\theta}\rightarrow\Sigma_{\theta}$, which
was described in section 5. Notice that this map can be extended
to $\Lambda(\Gamma)$ as the identity.\\

We construct the {\it $q$-fold cyclic covering},
${\cal{C}}^{q}(\Lambda)$, {\it branched over the wild knot}
$\Lambda(\Gamma)$, as follows (see \cite{rolfsen}). We consider
$\overline{\Sigma_{\theta}}\times [0,1]$ and we identify $\partial
\overline{\Sigma_{\theta}}\times \{t\}$, $t\in (0,1]$ with $\partial
\overline{\Sigma_{\theta}}\times \{0\}$ via the identity map and
$\Sigma_{\theta}\times\{0\}$ with $\Sigma_{\theta}\times\{1\}$ via the
$q^{th}$-iterate $\psi^{q}$, of the characteristic map. The resulting
space is ${\cal{C}}^{q}(\Lambda)$.\\

The covering map $p:{\cal{C}}^{q}(\Lambda)\rightarrow\mathbb{S}^{3}$
is defined in the usual way (see \cite{rolfsen}).

\begin{rems}
\begin{enumerate}
\item The branched covering
  $p:{\cal{C}}^{q}(\Lambda)\rightarrow\mathbb{S}^{3}$ is regular.
\item The space ${\cal{C}}^{q}(\Lambda)$ is compact and  path-connected. 
\end{enumerate}
\end{rems}
\begin{prop}
The fundamental group of ${\cal{C}}^{q}(\Lambda)$ is infinitely generated.
\end{prop}

\hspace{-.67cm}{\it Proof.}
Let $K$ be a non-trivial tame fibered knot with fiber $S$. Suppose
that the fundamental group of the surface $S$ is
$$
\Pi_{1}(S)=\{a_{1},a_{2},\ldots,a_{n}:r_{1},\ldots,r_{l}\}
$$
Let $T$ be the pearl necklace subordinate to $K$ and let $\psi$ be the
characteristic map of $\Lambda$. Then
$$
\hspace{-2cm}\Pi_{1}(\mathbb{S}^{3}-|T|)\cong\Pi_{1}(K)\cong\Pi_{1}(\mathbb{S}^{1})
\ltimes_{\psi_{\#}|_{T}}\Pi_{1}(S)
$$
$$
\hspace{3.2cm}=\{c,a_{1},a_{2},\ldots,a_{n}:r_{1},\ldots,r_{l},a_{k}*c=\psi_{\#}|_{T}(a_{k})\},
$$
and the fundamental group of the $q$-fold cyclic branched covering
${\cal{C}}^{q}(K)$ of $\mathbb{S}^{3}$ with branched set $K$ is 
$$
\Pi_{1}({\cal{C}}^{q}(K))=\{c,a_{1},a_{2},\ldots,a_{n}:r_{1},\ldots,r_{l},a_{k}*c=\psi_{\#}|_{T}(a_{k}),
\psi_{\#}^{q}|_{T}=1\}.
$$
By Van-Kampen's theorem, we have that the fundamental group of the
fiber $\Sigma_{\theta}$ is
$$
\Pi_{1}(\Sigma_{\theta})=\{a_{i1},a_{i2},\ldots,a_{in}:r_{i1},\ldots,r_{il}\}
$$
where $i\in\mathbb{N}$. Thus, the fundamental group of
$\mathbb{S}^{3}-\Lambda$ is
$$ 
\hspace{-2.7cm}\Pi_{1}(\mathbb{S}^{3}-\Lambda)\cong\Pi_{1}(\Lambda)\cong\Pi_{1}(\mathbb{S}^{1})
\ltimes_{\psi_{\#}}\Pi_{1}(\Sigma_{\theta})
$$
$$
\hspace{3.2cm}=\{c,a_{i1},a_{i2},\ldots,a_{in}:r_{i1},\ldots,r_{il},a_{ik}*c=\psi_{\#}(a_{ik})\},
$$
where $\psi_{\#}(a_{ik})$ is a product of
$a_{i1},a_{i2},\ldots,a_{in}$. This implies that the fundamental group of ${\cal{C}}^{q}(\Lambda)$ is
$$
\Pi_{1}({\cal{C}}^{q}(\Lambda))=\{c,a_{i1},a_{i2},\ldots,a_{in}:r_{i1},\ldots,r_{il},
a_{ik}*c=\psi_{\#}(a_{ik}),\psi_{\#}^{q}=1\}
$$
$$
\cong(\cdots((\Pi_{1}({\cal{C}}^{q}(K))*_{\{c\}}\Pi_{1}({\cal{C}}^{q}(K)))
*_{\{c\}}\cdots *_{\{c\}}\Pi_{1}({\cal{C}}^{q}(K)))
*_{\{c\}} \cdots
$$
which is an infinite free product of the fundamental group of
${\cal{C}}^{q}(K)$. Hence $\Pi_{1}({\cal{C}}^{q}(\Lambda))$ is
infinitely generated. $\blacksquare$

\begin{ex}
Let $K=T_{2,3}$ be the right-handed trefoil knot. Then its monodromy is order
six, up to isotopy. Consider the $5^{th}$-fold cyclic branched
covering ${\cal{C}}^{5}(\Lambda)$ of $\mathbb{S}^{3}$ with branched
set the wild knot $\Lambda$, obtained from $K$. Then
$$
\Pi_{1}({\cal{C}}^{5}(\Lambda))\cong(\cdots((\Pi_{1}({\cal{C}}^{5}(K))*_{\{c\}}\Pi_{1}({\cal{C}}^{5}(K)))
*_{\{c\}}\cdots *_{\{c\}}\Pi_{1}({\cal{C}}^{5}(K)))
*_{\{c\}} \cdots
$$
is the infinite free product of the fundamental group of
${\cal{C}}^{5}(K)$. In other words, it is the infinite free product of the binary icosahedral group.\\

On the other hand, by the Hurewicz homomorphism, we have than the first homology group of  
${\cal{C}}^{5}(\Lambda)$ is
$$
H_{1}({\cal{C}}^{5}(\Lambda),\mathbb{Z})\cong 
\Pi_{1}({\cal{C}}^{5}(\Lambda))/[\Pi_{1}({\cal{C}}^{5}(\Lambda)),\Pi_{1}({\cal{C}}^{5}(\Lambda))].
$$
We recall that ${\cal{C}}^{5}(K)$ is the Poincar\'e Sphere $\Sigma_{2,3,5}$, hence
$H_{1}({\cal{C}}^{5}(K),\mathbb{Z})=0.$
This implies that $H_{1}({\cal{C}}^{5}(\Lambda),\mathbb{Z})=0$.
\end{ex}
\begin{theo}
The space ${\cal{C}}^{q}(\Lambda)$ is not semilocally simply
connected. In particular, it does not admit a universal covering.
\end{theo}

\hspace{-.67cm}{\it Proof.}
Let $K$ be a non-trivial tame fibered knot with fiber $S$, and let
$\Lambda$ be the limit set obtained from $K$ through the reflecting
process.\\

Let $x\in\Lambda\subset{\cal{C}}^{q}(\Lambda)$. Consider $U$ an open,
connected neighbourhood of $x$. Then we can think $U$ as a
$[0,1]$-family of connected neighbourhoods 
$U_{t}\subset\overline{\Sigma_{\theta}}\times\{t\}$
glued together along their boundaries, with $U_{0}$ identified with
$U_{1}$ via the $q^{th}$-iterate, $\psi^{q}$ of the characteristic map of $\Lambda$. \\

Let $U_{x}$ be the largest
neighbourhood of $x$ in $\overline{\Sigma_{\theta}}$ that satisfies 
$U_{x}\subset U_{t}$, $t\in\mathbb{S}^{1}$. Then we have that a
infinite number of copies of the surface $S$ is contained in
$Int(U_{x})$. Consider $\tilde{U}=U_{x}\times[0,1]/(y,1)\sim (\psi^{q}(y),0)\subset U$. Then
$\Pi_{1}(\tilde{U})$ is a subgroup of $\Pi_{1}({U})$ and $\tilde{U}$ contains an infinite number
of copies of ${\cal{C}}^{q}(K)$. This implies that (see the proof of the above theorem) 
$$
\Pi_{1}(U)\geq \{c,a_{i1},a_{i2},\ldots,a_{in}:r_{i1},\ldots,r_{il},
a_{ik}*c=\psi_{\#}(a_{ik}),\psi_{\#}^{q}=1\}
$$ 
$$
\cong(\cdots((\Pi_{1}({\cal{C}}^{q}(K))*_{\{c\}}\Pi_{1}({\cal{C}}^{q}(K)))
*_{\{c\}}\cdots *_{\{c\}}\Pi_{1}({\cal{C}}^{q}(K)))
*_{\{c\}} \cdots
$$
where $i\in\mathbb{N}$. Therefore, 
${\cal{C}}^{q}(\Lambda)$ is not semilocally simply
connected. $\blacksquare$
 
\begin{coro}
The space ${\cal{C}}^{q}(\Lambda)$ is not a manifold.
\end{coro}

\begin{rems}
\begin{enumerate}
\item The space ${\cal{C}}^{q}(\Lambda)-\Lambda$ is a 3-manifold.
\item The group of covering transformations of ${\cal{C}}^{q}(\Lambda)$ is
$\mathbb{Z}/q\mathbb{Z}$. Hence the quotient space 
${\cal{C}}^{q}(\Lambda)/\mathbb{Z}/q\mathbb{Z}\cong\mathbb{S}^{3}$
is a manifold.
\end{enumerate}
\end{rems}

\begin{theo}
Every closed, connected, orientable 3-manifold contains either a {\rm fibered}
wild knot or a {\rm fibered} wild link. 
\end{theo}

\hspace{-.67cm}{\it Proof.}
Let $K$ be the figure-eight knot. Let $T$ be a pearl-necklace
subordinate to $K$. Then the wild knot $\Lambda(\Gamma)$ obtained from
$K$ through the reflecting process is contained in $|T|$.\\ 

By \cite{montesinos}, we have that every closed, connected, orientable 
3-manifold $M$ is a n-fold branched covering of
$\mathbb{S}^{3}$ with branched set $K$. Let $p$ be the covering
map. Then $p^{-1}(|T|)\subset M$ is a finite number of copies of $|T|$ and in 
the interior of each copy of $|T|$, we have a copy of $\Lambda(\Gamma)$.\\

Therefore, $M$ contains a {\it fibered wild knot or link}.
$\blacksquare$ 

\begin{rem}
Even though any closed, connected, orientable 3-manifold contains a
wild knot it is not so clear
that this wild knot can be chosen to be fibered.
\end{rem}

\end{document}